\documentclass{article}
\usepackage{amssymb,amsmath,amsthm,hyperref}
\usepackage[title]{appendix}

\usepackage{xcolor}
\usepackage{wasysym}

\newcommand\doublecheck{\textcolor{blue}{\checked\kern-0.6em\checked}}
\newcommand{\newvtheorem}[2]{\newtheorem{#1}[theorem]{\llap{\textnormal{\doublecheck} }#2}}
\newcommand{\newvtheoremstar}[2]{\newtheorem*{#1}{\llap{\textnormal{\doublecheck} }#2}}
\newvtheorem{vtheorem}{Theorem}
\newvtheorem{vlemma}{Lemma}
\newvtheoremstar{vremark}{Remark}

\begin{document}

\vspace*{-2cm}

\Large
 \begin{center}
The cardinality of a set containing the pairwise sums of a fixed number of integers \\ 

\hspace{10pt}

\large
Wouter van Doorn \\

\hspace{10pt}

\end{center}

\hspace{10pt}

\normalsize

\vspace{-25pt}
\centerline{\bf Abstract}
Revisiting a $50$-year-old estimate of Choi, Erd\H{o}s and Szemer\'edi, we show that if $A \subseteq \{1, 2, \ldots, 2n\}$ satisfies $|A| \ge n + 1.2 \cdot 10^8$, then there exist five distinct integers whose pairwise sums are all contained in $A$. In order to guarantee pairwise sums of three or four integers instead, we show that one can replace the constant $1.2 \cdot 10^8$ by $1$ or $3$ respectively, which are both optimal.

\section{Introduction}
For positive integers $n$ and $k \ge 3$, let $g_k(n)$ be the smallest integer such that for all sets $A \subseteq \{1, 2, \ldots, 2n\}$ with $|A| \ge n + g_k(n)$, one can find distinct integers $b_1, b_2, \ldots, b_k$ (not necessarily contained in $A$) with $b_i + b_j \in A$ for $1 \le i < j \le k$. Choi, Erd\H{o}s and Szemer\'edi introduced this function in \cite{ces}, and proved many bounds on it (which were already announced three years prior in \cite[p. 83]{Er72}). For example, they proved 
\begin{align*}
g_3(n) &\le 2, \\
g_4(n) &= O(1), \\
g_5(n) &= O\big(\log(n)\big),\\
g_6(n) &\asymp \sqrt{n}.
\end{align*}

They did not necessarily attempt to optimize all the implied constants however, and estimating the value of $g_k(n)$ is now listed as Erd\H{o}s Problem \#866 at \cite{bloom}. In this paper we revisit the definition of $g_k$, and prove that for all $n \ge 3$ we have $g_3(n) = 1$, $g_4(n) = 3$ and $g_5(n) < 1.2 \cdot 10^8$.

\section{Declaration of AI usage} \label{ai}
In this research ChatGPT (model: GPT-5.3 Instant) was initially used for brainstorming ideas. Most notably, it autonomously came up with the proof of Theorem \ref{g4} that we outline below. We further used the automated theorem proving tool Aristotle from Harmonic \cite{ari} to obtain Lean formalizations of essentially all results in this paper. More precisely, every stated result with a blue checkmark symbol has been successfully formalized, and the Lean file containing these formal proofs can be found at the author's GitHub page \cite{wvd}. Particularly noteworthy here is that, in the process of formalizing the result that $h_4(n) \le 3166$ for all $n$,\footnote{See Section \ref{def} for a definition of $h$, and Section \ref{finalthoughts} for a proof sketch of the claimed upper bound.} Aristotle managed to come up with further improvements, eventually showing $h_4(n) \le 2270$.

\section{Clarifying the main definition} \label{def}
Regarding the definition of $g_k(n)$, in all theorem statements in \cite{ces} the elements in $A \subseteq \{1, 2, \ldots, 2n\}$ are referred to as \emph{positive} integers, whereas the $b_i$ (whose pairwise sums need to be contained in $A$) are just referred to as integers. It is never acknowledged however, that with this definition it is not actually required that the $b_i$ are positive themselves as well. Of course, if both $b_1 \le 0$ and $b_2 \le 0$ (say), then $b_1 + b_2 \le 0 \notin A$, from which we deduce that at most one of the $b_i$ can be non-positive. But allowing (at most one of) the $b_i$ to be non-positive does actually matter! For example, it is claimed that the set 

\begin{equation} \label{adef}
A := \{1, 3, \ldots, 2n-1\} \cup \{2, 4, \ldots, 2^{\left \lfloor \log_2 2n \right \rfloor}\} 
\end{equation}
is an example for which

\begin{quote}
\emph{one cannot choose $b_1, \ldots, b_5$ such that $b_i + b_j$ ($1 \le i < j \le 5$) are all in $A$.} \cite[p. 40]{ces}
\end{quote}

In particular, this would give the lower bound $g_5(n) > \log_2 n$. This is incorrect as stated, however. With $$b_1 := -1, b_2 := 2, b_3 := 3, b_4 := 5, b_5 := 6,$$ one can verify that all pairwise sums do belong to $A$ whenever $n \ge 6$, contradicting the above quote. \\

We are therefore led to define the function $h_k(n)$ as the smallest integer such that for all sets $A \subseteq \{1, 2, \ldots, 2n\}$ with $|A| \ge n + h_k(n)$, one can find $k$ distinct \emph{positive} integers whose $\binom{k}{2}$ pairwise sums are contained in $A$. With this definition, later on we will see that we do recover the lower bound $h_5(n) > \log_2 n$, but that $g_k(n)$ and $h_k(n)$ are in general not the same function. That being said, the fact that (even if we only assume $b_i \in \mathbb{Z}$) at most one of the $b_i$ is non-positive, does imply the bounds 

\begin{equation} \label{generalgvsh}
g_k(n) \le h_k(n) \le g_{k+1}(n)
\end{equation}

for all $k \ge 3$ and all $n \ge 1$.

\section{Sums of three or four integers}
As we mentioned, the proposed lower bound example for $g_5(n)$ does not work if we do not put further restrictions on the $b_i$. Similarly, the example $$A := \{1, 3, \ldots, 2n-1\} \cup \{2 \}$$ given in \cite{ces} only helps to show $h_3(n) = 2$ for all $n \ge 4$. For $g_3(n)$ on the other hand, we actually claim that it is equal to $1$.

\begin{vtheorem} \label{g3} 
We have $g_3(n) = 1$ for all $n \ge 3$.
\end{vtheorem}

\begin{proof}
The lower bound $g_3(n) \ge 1$ follows from taking $A$ to be the set of all odd integers in $\{1, 2, \ldots, 2n\}$, as at least one sum $b_i + b_j$ must be even by the pigeonhole principle. \\

For the corresponding upper bound, let $A \subseteq \{1, 2, \ldots, 2n\}$ be any set with $|A| \ge n+1$, and assume $n \ge 3$. Then, once again by the pigeonhole principle, $A$ contains $\{m, m+1\}$ for some odd $m \in \{1, 3, \ldots, 2n-1\}$. If $A \setminus \{m \}$ only contains even integers, then $$\{2, 4, 6\} \subseteq \{2, 4, \ldots, 2n\} \subset A,$$ and we can choose $$b_1 := 0, b_2 := 2, b_3 := 4.$$ On the other hand, if $A$ contains an odd integer $2j+1 \neq m$, set $$b_1 := j, b_2 := j+1, b_3 := m-j.$$ We then need to show that the $b_i$ are distinct and that their pairwise sums are contained in $A$. First of all, it is clear that $b_1 \neq b_2$. Secondly, $b_1 \neq b_3$ as we would otherwise have $m = 2j$ contradicting the fact that $m$ is odd. Thirdly, $b_2 \neq b_3$ by the assumption $2j+1 \neq m$. And finally, the set of pairwise sums is $\{2j+1, m, m+1\}$ which is indeed a subset of $A$. 
\end{proof}

\begin{vremark} \label{g3small}
The examples $A := \{1, 2\}$ and $A := \{2, 3, 4\}$ show $g_3(1) = g_3(2) = 2$.
\end{vremark}

As it turns out (and despite our insistence on allowing the $b_i$ to be negative), for Theorem \ref{g3} we do not actually need to use any negative integers.

\begin{vtheorem} \label{g3nonnegative}
If $n \ge 3$ and $A \subseteq \{1, 2, \ldots, 2n\}$ is a set with $|A| \ge n+1$, then integers $0 \le b_1 < b_2 < b_3$ exist whose pairwise sums are contained in $A$.
\end{vtheorem}

\begin{proof}
We will prove it by induction; % we are even allowed to assume ${1, 2, 2n-1, 2n} \subset A$ by induction, but don't fully make use of this 
the interested reader is invited to verify it for $n = 3$. We then let $$A = \{a_1, a_2, \ldots, a_r\} \subseteq \{1, 2, \ldots, 2n\}$$ be any set with at least $n+1$ elements, and assume that the theorem is proven for all sets $A' \subseteq \{1, 2, \ldots, 2n-2\}$ with at least $n$ elements. Now, if $1 \notin A$, define $$A' := \{a_1-2, a_2-2, \ldots, a_r-2 \} \setminus \{0 \}.$$ Then the induction hypothesis applies to $A'$, so let $b'_1, b'_2, b'_3$ be three non-negative integers whose pairwise sums are in $A'$. Then we can simply choose $b_i := b'_i + 1$, as $b'_i + b'_j \in A'$ implies $b_i + b_j = b'_i + b'_j + 2 \in A$. On the other hand, if $1 \in A$, then there are two options. Either $\{m, m+1\} \subset A$ for some $m \ge 2$, in which case we can set $b_1 := 0, b_2 := 1, b_3 := m$. Alternatively, if $A$ does not contain any pair of consecutive integers besides possibly $\{1, 2\}$, then $A$ must be exactly equal to $\{1, 2, 4, 6, \ldots, 2n\}$. In that case we can set $b_1 := 0, b_2 := 2, b_3 := 4$.
\end{proof} 

From $k = 3$ we now move on to $k = 4$.

\begin{vtheorem} \label{g4}
We have $g_4(n) = 3$ for all $n \ge 2$.
\end{vtheorem}

\begin{proof}
For the lower bound, we choose $$A := \{1, 3, \ldots, 2n-1\} \cup \{2n-2, 2n \}.$$ To see why $A$ does not contain the pairwise sums of any four integers $b_1, b_2, b_3, b_4$, we are free to assume that $b_1$ and $b_2$ are even, and $b_3$ and $b_4$ are both odd. Indeed, if three of them have the same parity, then their three pairwise sums are all even, while $A$ only contains two even integers. Furthermore assume without loss of generality $b_1 < b_2$ and $b_3 < b_4$. Now, because the smallest even integer in $A$ is $2n-2$, we must have that both $b_1 + b_2$ and $b_3 + b_4$ are at least $2n-2$. In particular, $b_2$ and $b_4$ are both at least $n$, which means their sum is at least $2n+1$ (as they cannot both be equal to $n$), which is not an element of $A$. \\

As for the upper bound, let $A \subseteq \{1, 2, \ldots, 2n\}$ be any set with $|A| \ge n+3$. Since the result is vacuous for $n = 2$, we may assume $n \ge 3$. The pigeonhole principle then implies that there exist three distinct even integers $a_1, a_2, a_3$ such that $A$ contains $\{a_i, 2n+1-a_i\}$ for all $1 \le i \le 3$. We now set 
\begin{align*}
b_1 &:= \frac{a_1 + a_2 - a_3}{2}, \\
b_2 &:= \frac{a_1 + a_3 - a_2}{2}, \\
b_3 &:= \frac{a_2 + a_3 - a_1}{2}, \\
b_4 &:= \frac{4n+2 - a_1 - a_2 - a_3}{2} .
\end{align*}

These integers are defined so that we have
\begin{align*}
b_1 + b_2 &= a_1, \\
b_1 + b_3 &= a_2, \\
b_1 + b_4 &= 2n+1-a_3, \\
b_2 + b_3 &= a_3, \\
b_2 + b_4 &= 2n+1-a_2, \\
b_3 + b_4 &= 2n+1-a_1.
\end{align*}

In particular, $b_i + b_j \in A$ for all $1 \le i < j \le 4$. Finally, $b_1, b_2, b_3$ are all distinct by the assumption that $a_1, a_2, a_3$ are distinct as well, while $b_4 \notin \{b_1, b_2, b_3\}$ by the fact that $b_4$ has the opposite parity.
\end{proof}

\section{Lower bounds for sums of five integers}
Let us first formally verify that the example Choi, Erd\H{o}s and Szemer\'edi propose, does work to provide a lower bound on $h_5(n)$.

\begin{vtheorem}[Choi-Erd\H{o}s-Szemer\'edi] \label{h5lower}
We have $h_5(n) > \log_2 n$ for all $n \in \mathbb{N}$.
\end{vtheorem}

\begin{proof}
With $A$ as in \eqref{adef}, we note that we indeed have $|A| > n + \log_2 n$ for all $n$. Hence, with $b_1, b_2, \ldots, b_5$ five arbitrary distinct positive integers, it is sufficient to show that at least one of the sums $b_i + b_j$ (with $i,j$ distinct) is not contained in $A$. Now, by the pigeonhole principle at least three of the $b_i$ have the same parity. \\

Let us first assume that there are three $b_i$ that are all even, so that without loss of generality we may write $$b_1 = 2^{l_1} m_1, b_2 = 2^{l_2} m_2, b_3 = 2^{l_3} m_3$$ with $m_1$, $m_2$, and $m_3$ odd, and $1 \le l_1 \le l_2 \le l_3$. Now, if one of the latter two inequalities is strict, then $$b_1 + b_3 = 2^{l_1} (m_1 + 2^{l_3 - l_1} m_3)$$ is an even integer which is not a power of two, and hence not in $A$. If $l_1 = l_2 = l_3$ on the other hand, then two of the $m_i$ (say $m_1$ and $m_2$) are congruent modulo $4$, in which case $$b_1 + b_2 = 2^{l_1} (m_1 + m_2)$$ is an even integer which is not a power of two, as $m_1 + m_2 \equiv 2 \pmod{4}$ with $m_1 + m_2 > 2$. \\

Since the latter argument still works when $l_1 = l_2 = l_3 = 0$, this also deals with the case where at least three of the $b_i$ are odd, finishing the proof.
\end{proof}

Let us now return to $g_k$, for which we have the following, seemingly modest, lower bound for $k = 5$.

\begin{vtheorem} \label{g5lower}
We have $g_5(n) \ge 4$ for all $n \ge 3$. 
\end{vtheorem}

\begin{proof}
More precisely, we will prove that for the set 

\begin{equation} \label{adef2}
A := \{1, 3, \ldots, 2n-1\} \cup \{2n-4, 2n-2, 2n\}, 
\end{equation}

there are no five distinct integers whose pairwise sums are all in $A$. So let $b_1, b_2, b_3, b_4, b_5$ be any distinct integers and assume by contradiction that $A$ does contain all pairwise sums. If at least four of them have the same parity, then there are at least five distinct even sums, while $A$ only contains three even integers. We may therefore assume without loss of generality $$b_1 \equiv b_2 \equiv b_3 \not \equiv b_4 \equiv b_5 \pmod{2}$$ with $$b_1 < b_2 < b_3 \qquad \text{and} \qquad b_4 < b_5.$$ Given the fact that $2n-4, 2n-2$ and $2n$ are the only even integers in $A$, we then deduce 
\begin{align*}
b_1 + b_2 &= 2n-4, \\
b_1 + b_3 &= 2n-2, \\
b_2 + b_3 &= 2n, \\
b_4 + b_5 &\ge 2n-4.
\end{align*}

Solving the first three equations gives $$b_1 = n-3, b_2 = n-1, b_3 = n+1,$$ while the final inequality implies $b_5 \ge n-1$, so that we actually get $b_5 \ge n$ (as we already have $b_2 = n-1$). But then $b_3 + b_5 \ge 2n+1 \notin A$.
\end{proof}

The example set as defined in \eqref{adef2} is optimal in the following narrow sense.

\begin{vtheorem} \label{g5special}
If $A$ is a set with $$\{1, 3, \ldots, 2n-1\} \subset A \subseteq \{1, 2, \ldots, 2n\}$$ and $|A| \ge n + 4$, then distinct integers $b_1, b_2, \ldots, b_5$ exist such that $b_i + b_j \in A$ for all $1 \le i < j \le 5$.
\end{vtheorem}

As neither the statement nor the proof of Theorem \ref{g5special} is particularly illuminating (especially in light of the main result in the next section), we have moved the proof to Appendix \ref{appendix}

\section{A uniform upper bound for sums of five integers}
Using arguments that are essentially borrowed from \cite{ces}, this section is aimed at proving that $g_5(n)$ is uniformly bounded. We first need a preliminary result that upper bounds $g_k(n)$ if we restrict to sets of even integers. To be able to state it, for $x \ge 1$ we define 

\begin{equation} \label{f3def}
f_3(x) := \sqrt{\frac{x}{2}} + \left(\frac{x}{2}\right)^{1/4} + \frac{1}{2}
\end{equation}

and, for $k \ge 4$, we further define 

\begin{equation} \label{fkdef}
f_{k}(x) := \sqrt{2xf_{k-1}(x) + \frac{1}{4}} + \frac{1}{2}.
\end{equation}

We then have the following slightly more optimized version of the corollary to Lemma $A$ in \cite{ces}.

\begin{vlemma} \label{ceslemgeneral}
If $A_0 = \{2 \le a_1 < a_2 < \cdots < a_r\} \subset \mathbb{N}$ is a set of even integers for which the inequality $r \ge f_k(a_r - a_1)$ holds for some $k \ge 3$, then distinct integers $b_1, b_2, \ldots, b_k$ exist such that $A_0$ contains all $\binom{k}{2}$ pairwise sums $b_i + b_j$.
\end{vlemma}

\begin{proof}
We will prove by induction that $r \ge f_k(a_r - a_1)$ implies the existence of integers $c_1, c_2, \ldots, c_k$ with $c_2, c_3, \ldots, c_k$ non-zero and pairwise distinct, such that $A_0$ contains all subset sums containing $c_1$. That is, for any $S \subseteq \{1, 2, \ldots, k\}$ with $1 \in S$, we have $\sum_{i \in S} c_i \in A_0$. If such $c_i$ indeed exist, then the lemma follows by taking $b_1 := \frac{1}{2}c_1$ and $b_i := \frac{1}{2}c_1 + c_i$ for $2 \le i \le k$. \\

We first need to prove it for $k = 3$. In this case, if a positive integer $m$ and indices $i_1 < i_2$, $j_1 < j_2$ exist with $$a_{j_1} - a_{i_1} = a_{j_2} - a_{i_2} = m,$$ then one can choose $$c_1 := a_{i_2}, c_2 := a_{i_1} - a_{i_2}, c_3 := m,$$ where we have $c_2 < 0 < c_3$. However, if such an $m$ and indices do not exist, then $\{a_1/2, a_2/2, \ldots, a_r/2\}$ forms a so-called Sidon set, for which the bound $$r < \sqrt{\frac{a_r}{2} - \frac{a_1}{2}} + \left(\frac{a_r}{2} - \frac{a_1}{2}\right)^{1/4} + \frac{1}{2} = f_3(a_r - a_1)$$ was proved in \cite{sidon}. Now let $k > 3$ and assume that the lemma holds for $k-1$. \\

There are a total of $\frac{1}{2}r(r-1)$ differences $a_j - a_i$ in the interval $[2, a_r-a_1]$, and all of them are even. Hence, by the pigeonhole principle, if 

\begin{equation} \label{inductionineq}
\frac{1}{2}r(r-1) \ge (a_r-a_1)f_{k-1}(a_r-a_1),
\end{equation}

then there are positive integers $m$ and $s \ge f_{k-1}(a_r-a_1)$ such that $m$ can be written in $s$ disjoint ways as a difference $$a_{j_1} - a_{i_1} = a_{j_2} - a_{i_2} = \cdots = a_{j_{s}} - a_{i_{s}} = m.$$ And indeed, with definition \eqref{fkdef} in mind, inequality \eqref{inductionineq} is precisely equivalent to $r \ge f_k(a_r - a_1)$ by the quadratic formula. For simplicity, write $$S := \{a_{i_1}, a_{i_2}, \ldots, a_{i_{s}}\} \qquad \text{and} \qquad T := \{a_{j_1}, a_{j_2}, \ldots, a_{j_{s}}\}.$$ Now we apply the induction hypothesis to $S$, to obtain $c_1, c_2, \ldots, c_{k-1}$ for which all the subset sums containing $c_1$ belong to $S$. We then set $c_k := m$ and note that $c_k$ is distinct from the integers $c_2, c_3, \ldots, c_{k-1}$, because otherwise $c_1 + m \in S \cap T$, contradicting the fact that $S$ and $T$ are disjoint. With $c_k = m$ included, we deduce that all the subset sums containing $c_1$ now belong to $S \cup T \subseteq A_0,$ as desired.
\end{proof}

We are now ready to prove our main theorem, by noting that from the definitions in \eqref{f3def} and \eqref{fkdef}, one can verify by induction that for all $k \ge 3$ the inequalities\footnote{In fact, for the first inequality the explicit estimate of $$f_k(x) \le 2^{1 - \frac{3}{2^{k-2}}}x^{1 - \frac{1}{2^{k-2}}} + 2x^{1 - \frac{3}{2^{k-1}}}$$ is possible.}

\begin{equation} \label{fkupper}
f_k(x) \le \big(2 + o(1)\big)x^{1 - \frac{1}{2^{k-2}}}
\end{equation}

and 

\begin{equation} \label{fksublinear}
2f_k(x) > f_k(2x)
\end{equation}

hold for all $x \ge 1$.

\begin{vtheorem} \label{g5upper}
We have $g_5(n) < 1.2 \cdot 10^8$ for all $n \in \mathbb{N}$. 
\end{vtheorem}

\begin{proof}
Let $C$ be the smallest integer for which the inequality 

\begin{equation} \label{xineq}
\frac{x}{12} + \frac{C}{2} - 2 > f_5(x)
\end{equation}

holds for all $x \ge 1$. We remark that such a $C$ certainly exists, as the left-hand side grows linearly, while the right-hand side has sublinear growth by the upper bound \eqref{fkupper}. In fact, calculating the derivative of $\frac{x}{12} - f_5(x)$ and setting it equal to $0$ gives $$C = 113,591,719 < 1.2 \cdot 10^8.$$ Now let $A \subseteq \{1, 2, \ldots, 2n\}$ be any set with $|A| \ge n + C$, let $A_0 \subset A$ be the subset of all even integers in $A$, and write $|A_0| = t + C$ for some $t \ge 0$. If $t \ge \frac{n}{6}$, then $$|A_0| \ge \frac{n}{6} + C > f_5(2n)$$ for all $n \ge 1$ by multiplying inequality \eqref{xineq} by $2$ and applying \eqref{fksublinear}, which implies we are done by Lemma \ref{ceslemgeneral}. We may therefore assume $t < \frac{n}{6}$ for the rest of this proof. \\

Now, if there are at most $4$ even elements of $A$ contained in the closed interval $[6t+2, 2n-6t-2]$, then $t > 0$ and there are either at least $$r := \left \lceil \frac{|A_0|-4}{2} \right \rceil = \left \lceil \frac{t}{2} + \frac{C}{2} - 2 \right \rceil$$ even elements of $A$ in the interval $[2, 6t]$, or at least $r$ even elements of $A$ in $[2n-6t, 2n]$. In either situation, with $x := 6t$ we have an interval of length at most $x$ with at least $r$ even elements. Lemma \ref{ceslemgeneral} then finishes the proof again, as $$r \ge \frac{t}{2} + \frac{C}{2} - 2 = \frac{x}{12} + \frac{C}{2} - 2 > f_5(x),$$ by inequality \eqref{xineq}. \\

Hence, without loss of generality we may assume that there are at least $5$ even elements of $A$ contained in the interval $[6t+2, 2n-6t-2]$. In particular, we may assume the existence of three even integers $a_1, a_2, a_3$ with $$6t+2 \le a_1 < a_2 < a_3 \le n$$ or three even integers $a'_1, a'_2, a'_3$ with $$n \le a'_1 < a'_2 < a'_3 \le 2n - 6t - 2.$$ Similar to the proof of Theorem \ref{g4}, we then define
\begin{align*}
b_1 &:= \frac{a_1 + a_2 - a_3}{2}, \\
b_2 &:= \frac{a_1 + a_3 - a_2}{2}, \\
b_3 &:= \frac{a_2 + a_3 - a_1}{2}
\end{align*}

in the first case, or
\begin{align*}
b'_1 &:= \frac{a'_1 + a'_2 - a'_3}{2}, \\
b'_2 &:= \frac{a'_1 + a'_3 - a'_2}{2}, \\
b'_3 &:= \frac{a'_2 + a'_3 - a'_1}{2}
\end{align*}

in the second case. As before, 
\begin{align*}
b_1 &< b_2 < b_3, \\
b_1 &\equiv b_2 \equiv b_3 \pmod{2}, \\
\{b_1& + b_2, b_1 + b_3, b_2 + b_3\} \subset A.
\end{align*}

And similarly for $b'_1, b'_2, b'_3$. \\

In the first case we consider the set of all integers $p, q$ with $b_1 \not \equiv p \equiv q \pmod{2}$ such that $p+q = a_3$ and $$\frac{a_3}{2} - 6t - 2 \le p < q \le \frac{a_3}{2} + 6t+2,$$ and in the second case we consider the set of all integers $p, q$ with $b'_1 \not \equiv p \equiv q \pmod{2}$ such that $p+q = a'_1$ and $$\frac{a'_1}{2} - 6t - 2 \le p < q \le \frac{a'_1}{2} + 6t+2.$$ In both cases there are at least $3t+1$ such pairs $(p,q)$, while we claim that all pairwise sums are contained in $[1, 2n]$. Indeed $p+q \in \{a_3, a'_1\} \subset [1, 2n]$, while
\begin{align*}
p + b_1 &\ge \frac{a_3}{2} - 6t - 2 + \frac{a_1 + a_2 - a_3}{2} \\
&\ge - 6t - 2 + \frac{(6t+2) + (6t+4)}{2} \\
&= 1, \\
q + b_3 &\le \frac{a_3}{2} + 6t + 2 + \frac{a_2 + a_3 - a_1}{2} \\
&\le \frac{n}{2} + 6t + 2 + \frac{(n-2) + n - (6t+2)}{2} \\
&= \frac{3n}{2} + 3t \\
&< 2n
\end{align*}

in the first case, and 
\begin{align*}
p + b'_1 &\ge \frac{a'_1}{2} - 6t - 2 + \frac{a'_1 + a'_2 - a'_3}{2} \\
&\ge \frac{n}{2} - 6t - 2 + \frac{n + (n+2) - (2n - 6t-2)}{2} \\
&= \frac{n}{2} - 3t \\
&> 0, \\
q + b'_3 &\le \frac{a'_1}{2} + 6t + 2 + \frac{a'_2 + a'_3 - a'_1}{2} \\
&\le 6t + 2 + \frac{(2n - 6t - 4) + (2n - 6t - 2)}{2} \\
&= 2n-1
\end{align*}
in the second case. \\

Now, $A$ misses at most $|A_0| - C = t$ odd integers in the interval $[1, 2n]$, while in either case we have at least $3t+1$ pairs $(p, q)$ to work with. Hence, in the first case there must be a pair $(p, q)$ such that $$b_1 + p, b_2 + p, b_3 + p, b_1 + q, b_2 + q, b_3 + q$$ are all in $A$, as each missing odd integer rules out at most three pairs. We may therefore take $b_4 := p, b_5 := q$. The analogous argument replacing $b_1, b_2, b_3$ by their primed counterparts works for the second case. 
\end{proof}

\section{Concluding thoughts and remarks} \label{finalthoughts}
In \cite{ces} it was proved that $h_4(n) = O(1)$. For an explicit variant, note that Theorem \ref{g5upper} implies by inequality \eqref{generalgvsh} that we actually have $h_4(n) < 1.2 \cdot 10^8$. However, a much better bound is possible. For this, we recall that Ruzsa \cite[Theorem 4.6]{imre} proved that for any set $A \subseteq \{1, 2, \ldots, n\}$ with $$|A| > \sqrt{n} + 4n^{1/4} + 11$$ there are four distinct integers $a_1 < a_2 < a_3 < a_4$ in $A$ for which $a_1 + a_4 = a_2 + a_3$. Sets for which no such $a_i$ exist are called weak Sidon sets, and note the difference with the definition of a Sidon set, where even $a_2 = a_3$ is not allowed. In any case, with four such integers, we can define $$c_1 := a_1, c_2 := a_2 - a_1, c_3 := a_3 - a_1$$ in the proof of the base case of Lemma \ref{ceslemgeneral}, to obtain positive $c_i$ and hence positive $b_i$. With the definition in \eqref{f3def} appropriately altered to $$F_3(x) := \sqrt{\frac{x}{2}} + 4\left(\frac{x}{2}\right)^{1/4} + 11$$ and $F_4(x)$ analogously defined according to \eqref{fkdef}, the rest of the proof of Lemma \ref{ceslemgeneral} then still goes through the exact same way. \\

To use this to find an upper bound for $h_4(n)$, we sketch an argument similar to the proof of Theorem \ref{g5upper}; let $A \subseteq \{1, 2, \ldots, 2n\}$ be a set with $|A| \ge n + C$ for some constant $C$ that we have yet to determine, let $A_0 \subset A$ be the subset of all even integers in $A$, and write $|A_0| = t + C$ for some $t \ge 0$. If there exists an even element $$2m \in [8t+6, 2n-8t-6] \cap A_0,$$ then first choose $b_1 := m-1$, $b_2 := m+1$. Secondly, choose $b_3 \in [m-4t-2, m-2]$ and $b_4 := 2m - b_3$ with $b_3 \equiv b_4 \not \equiv b_1 \pmod{2}$ and such that $$b_1 + b_3, b_1 + b_4, b_2 + b_3, b_2 + b_4$$ are all elements of $A$. The latter is possible as there are at least $2t+1$ possible choices for $b_3$, while every odd integer not contained in $A$ eliminates at most two possibilities, and $A$ misses at most $t$ odd integers. \\

On the other hand, if $[8t+6, 2n-8t-6] \cap A_0$ is empty, then we have at least $\left \lceil \frac{t}{2} + \frac{C}{2} \right \rceil$ even elements in an interval of length $x := 8t+6$. Hence, choosing $C$ large enough so that $$\frac{x-6}{16} + \frac{C}{2} > F_4(x)$$ for all $x \ge 1$ would finish the proof, and one can verify that $C = 3166$ does the job. We thereby find $h_4(n) \le 3166$ for all $n$, although the formalization that is available at \cite{wvd} actually provides a slightly better bound, as explained in Section \ref{ai}. Further improvements are plausible by e.g. applying stronger bounds on weak Sidon sets. \\

Turning to $k = 5$, despite the rather large constant in Theorem \ref{g5upper}, it could be the case that $g_5(n) \le 5$ holds for all $n$. We have tried to look for counterexamples up to $n = 15$, but none were found. In fact, we cannot even exclude the possibility that $g_5(n) \le 4$ for all large enough $n$, which would be optimal by Theorem \ref{g5lower}. We note however that e.g. $A = \{1, 2, 4, 5, 6, \ldots, 10\}$ does serve as an example that $g_5(n) \ge 5$ is at least possible for some small $n$. \\

As for $k = 6$, let $\{a_1, a_2, \ldots, a_r\} \subseteq \left\{1, 2, \ldots, \left \lfloor n/2 \right \rfloor\right\}$ be a Sidon set, and recall that such sets exist with $r = \left(\frac{1}{\sqrt{2}} + o(1)\right)\sqrt{n}$ (see e.g. \cite{obryant}). Now define $$A := \{1, 3, \ldots, 2n-1\} \cup \{4a_1 - 2, 4a_2 - 2, \ldots, 4a_r-2\}.$$ In \cite{ces} this set $A$ is used to prove the lower bound $$h_6(n) \ge \left(\frac{1}{\sqrt{2}} - o(1)\right)\sqrt{n},$$ but it works just as well as a lower bound for $g_6(n)$. As they also prove an upper bound with the same order of growth, we have
\begin{equation} \label{particulargvsh}
h_k(n) = O\big(g_k(n)\big)
\end{equation}

for $k \in \{3, 4, 6\}$, which one can compare with inequality \eqref{generalgvsh}. However, by Theorem \ref{h5lower} and Theorem \ref{g5upper}, inequality \eqref{particulargvsh} does of course not hold for $k = 5$. \\

Finally, similar to how inequality \eqref{fkupper} can be shown, for larger values of $k$ one can combine the aforementioned bound on weak Sidon sets with the definition in \eqref{fkdef} to prove the following general upper bound, which is a slight improvement over \cite[Theorem 5]{ces}.

\begin{vtheorem} \label{generalupper}
For all $k \ge 3$ there exists an $N$ such that for all $n \ge N$ we have $$g_k(n) \le h_k(n) < 4n^{1 - \frac{1}{2^{k-2}}}.$$
\end{vtheorem}

\begin{appendices}
\renewcommand{\thesection}{\Alph{section}.}
\section{Proof of Theorem \ref{g5special}} \label{appendix}
In this appendix we record the proof of Theorem \ref{g5special}, for completeness' sake.

\begin{proof}[Proof of Theorem \ref{g5special}]
Let $2 \le a_1 < a_2 < a_3 < a_4 \le 2n$ be four even integers in $A$ (which exist by the assumption $|A| \ge n+4$), and first assume $a_2 + a_3 \le 2n$. We then define 
\begin{align*}
b_1 &:= \frac{a_1 + a_2 - a_3}{2}, \\
b_2 &:= \frac{a_1 + a_3 - a_2}{2}, \\
b_3 &:= \frac{a_2 + a_3 - a_1}{2}, \\
b_4 &:= 1 - b_1, \\
b_5 &:= a_3 - b_4.
\end{align*}

We need to check that these $b_i$ are distinct and that the $10$ pairwise sums are in $A$. To see why the $b_i$ are distinct, it is quickly verified that 

\begin{equation} \label{b1lb2lb3}
b_1 < b_2 < b_3
\end{equation}

and 

\begin{equation} \label{mod2}
b_1 \equiv b_2 \equiv b_3 \not\equiv b_4 \equiv b_5 \pmod{2}.
\end{equation}

Hence, to prove distinctness it suffices to show 

\begin{equation} \label{b4lb5}
b_4 < b_5.
\end{equation}

And indeed, $$b_5 - b_4 = a_3 - 2b_4 = a_3 + 2b_1 - 2 = a_1 + a_2 - 2 > 0.$$

To prove that $A$ contains all $10$ pairwise sums, one can quickly verify that the four sums $b_1 + b_2$, $b_1 + b_3$, $b_2 + b_3$ and $b_4 + b_5$ are all in $\{a_1, a_2, a_3\} \subset A$. Furthermore, equation \eqref{mod2} implies that the other six pairwise sums are odd. As $A$ contains by assumption all odd integers, by equations \eqref{b1lb2lb3} and \eqref{b4lb5} it is sufficient to show the inequalities 

\begin{equation} \label{necineqs}
1 \le b_1 + b_4 < b_3 + b_5 \le 2n-1.
\end{equation}

The first inequality is immediate by the definition of $b_4$, while for the second inequality we calculate $$b_3 + b_5 = b_3 + a_3 + b_1 - 1 = a_2 + a_3 - 1 \le 2n-1,$$ where the final inequality uses the assumption $a_2 + a_3 \le 2n$. \\

The second case we have to deal with is $a_2 + a_3 > 2n$ and $a_3 < 2n-2$, in which case we define
\begin{align*}
b_1 &:= \frac{a_2 + a_3 - a_4}{2}, \\
b_2 &:= \frac{a_2 + a_4 - a_3}{2}, \\
b_3 &:= \frac{a_3 + a_4 - a_2}{2}, \\
b_4 &:= a_2 + b_3 - (2n-1), \\
b_5 &:= (2n-1) - b_3.
\end{align*}

Analogously, we now have $$\{b_1 + b_2, b_1 + b_3, b_2 + b_3, b_4 + b_5\} = \{a_2, a_3, a_4\} \subset A,$$ while both equation \eqref{b1lb2lb3} and \eqref{mod2} are still quickly verified. It is therefore similarly sufficient to show inequalities \eqref{b4lb5} and \eqref{necineqs}. And indeed:
\begin{align*}
b_5 - b_4 &= (4n-2) - a_2 - 2b_3 = (4n-2) - a_3 - a_4 > 0, \\
b_1 + b_4 &= b_1 + a_2 + b_3 - (2n-1) = a_2 + a_3 - (2n-1) \ge 1, \\
b_3 + b_5 &= 2n-1.
\end{align*}

The final case is where $a_2 + a_3 > 2n$ and $a_3 = 2n-2$. We then define
\begin{align*}
b_1 &:= \frac{a_2 + a_3 - a_4}{2}, \\
b_2 &:= \frac{a_2 + a_4 - a_3}{2}, \\
b_3 &:= \frac{a_3 + a_4 - a_2}{2}, \\
b_4 &:= a_1 + b_3 - (2n-1), \\
b_5 &:= (2n-1) - b_3.
\end{align*}

Once again, verifying inequalities \eqref{b4lb5} and \eqref{necineqs} is sufficient:
\begin{align*}
b_5 - b_4 &= (4n-2) - a_1 - 2b_3 = (4n-2) - a_1 + a_2 - a_3 - a_4 = a_2 - a_1 > 0, \\
b_1 + b_4 &= b_1 + a_1 + b_3 - (2n-1) = a_1 + a_3 - (2n-1) = a_1 - 1 \ge 1, \\
b_3 + b_5 &= 2n-1. \qedhere
\end{align*}
\end{proof}
\end{appendices}

\end{document}